\documentclass[12pt]{amsart}
\usepackage{amssymb}

\newtheorem{thm}{Theorem}[section]

\newtheorem{cor}{Corollary}[section]
\newtheorem{conj}{Conjecture}[section]
\newtheorem{prop}{Proposition}[section]
\newtheorem{question}{Question}[section]

\theoremstyle{remark}

\theoremstyle{definition}
\newtheorem{defn}{Definition}[section]

\numberwithin{equation}{section}

\allowdisplaybreaks[4]

\begin{document}

\newcommand{\claimref}[1]{Claim \ref{#1}}
\newcommand{\thmref}[1]{Theorem \ref{#1}}
\newcommand{\propref}[1]{Proposition \ref{#1}}
\newcommand{\lemref}[1]{Lemma \ref{#1}}
\newcommand{\coref}[1]{Corollary \ref{#1}}
\newcommand{\remref}[1]{Remark \ref{#1}}
\newcommand{\conjref}[1]{Conjecture \ref{#1}}
\newcommand{\questionref}[1]{Question \ref{#1}}
\newcommand{\secref}[1]{Sec. \ref{#1}}
\newcommand{\ssecref}[1]{\ref{#1}}
\newcommand{\sssecref}[1]{\ref{#1}}

\def \d#1{\displaystyle{#1}}
\def \mult{\mathop{\mathrm{mult}}\nolimits}
\def \rank{\mathop{\mathrm{rank}}\nolimits}
\def \codim{\mathop{\mathrm{codim}}\nolimits}
\def \Ord{\mathop{\mathrm{Ord}}\nolimits}
\def \Var{\mathop{\mathrm{Var}}\nolimits}
\def \Ext{\mathop{\mathrm{Ext}}\nolimits}
\def \EQ{\Leftrightarrow}
\def \Pic{\mathop{\mathrm{Pic}}\nolimits}
\def \Spec{\mathop{\mathrm{Spec}}\nolimits}
\def \mapright#1{\smash{\mathop{\longrightarrow}\limits^{#1}}}
\def \mapleft#1{\smash{\mathop{\longleftarrow}\limits^{#1}}}
\def \mapdown#1{\Big\downarrow\rlap{$\vcenter{\hbox{$\scriptstyle#1$}}$}}
\def \smapdown#1{\downarrow\rlap{$\vcenter{\hbox{$\scriptstyle#1$}}$}}
\def \A{{\mathbb A}}
\def \I{{\mathcal I}}
\def \J{{\mathcal J}}
\def \CO{{\mathcal O}}
\def \C{{\mathcal C}}
\def \BC{{\mathbb C}}
\def \m{{\mathcal M}}
\def \H{{\mathcal H}}
\def \S{{\mathcal S}}
\def \Z{{\mathcal Z}}
\def \BZ{{\mathbb Z}}
\def \Y{{\mathcal Y}}
\def \T{{\mathcal T}}
\def \P{{\mathbb P}}
\def \G{{\mathbb G}}
\def \F{{\mathbb F}}
\def \BR{{\mathbb R}}
\def \Hom{{\mathrm{Hom}}}
\def \nilrad{{\mathrm{nilrad}}}
\def \Supp{{\mathrm{Supp}}}
\def \Ann{{\mathrm{Ann}}}
\def \closure#1{\overline{#1}}
\def \EQ{\Leftrightarrow}
\def \imply{\Rightarrow}
\def \isom{\cong}
\def \embed{\hookrightarrow}
\def \tensor{\otimes}
\def \wt#1{{\widetilde{#1}}}
\def \Jac{{\mathrm{Jac}}}

\title{On the Intersection of Two Plane Curves}

\author{Xi Chen}

\address{Department of Mathematics\\
South Hall, Room 6607\\
University of California\\
Santa Barbara, CA 93106}
\email{xichen@math.ucsb.edu}
\date{March 13, 2000}
\maketitle

\section{Introduction and Statement of Results}

The following question was raised and partially answered by Geng Xu in \cite{X}.

\begin{question}\label{q1}
Let $D$ be a general degree $d$ curve in $\P^2$. What is the minimal number
$i(d, m)$ of points in the set-theoretical intersection $C\cap D$ for any
degree $m$ irreducible curve $C$ (suppose that $C$ and $D$ meet properly)?
\end{question}

This problem is related to a conjecture of Kobayashi and Zaidenberg which
states that for a sufficiently general curve $D\subset \P^2$ of degree
$d\ge 5$, as general in the sense that $D$ lies in $|\CO_{\P^2}(d)|\isom
\P^{d(d+3)/2}$ with countably many closed proper subvarieties removed, the
affine variety $\P^2\backslash D$ is
hyperbolic. One necessary condition for $\P^2\backslash D$ being hyperbolic
is that there is no rational curve $C\subset \P^2$ meeting
$D$ set-theoretically at fewer than three points; otherwise, there is going
to be a nonconstant holomorphic map $\BC\to C\backslash (C\cap D) \subset
\P^2\backslash D$. This property of $\P^2\backslash D$ was called
``algebraic hyperbolic'' in \cite{DSW}.

\begin{defn}
A quasi-projective variety is called {\it algebraic hyperbolic\/} if it
does not contain a curve whose normalization is an
elliptic curve or a rational curve with two points removed, i.e.,
$\P^1\backslash \{p, q\}\isom\BC^*\isom\Spec \BC[x, x^{-1}]$.
\end{defn}

Obviously, hyperbolicity implies algebraic hyperbolicity for
smooth quasi-projective varieties.

Using an elegant deformation-theoretical argument, Xu proved the following
\cite[Theorem 1]{X}.

\begin{thm}[Xu]\label{txu}
For $d\ge 3$, $\min_{m > 0} i(d, m) = d-2$.
\end{thm}
 
He thus concluded that every curve $C\subset \P^2$ meets $D$ at
no less than three distinct points and hence $\P^2\backslash D$ is
algebraic hyperbolic for a sufficiently general curve
$D$ of degree $d\ge 5$. This bound is sharp for $m = 1$ and it is achieved
by a bitangent or flex line to $D$. 

The purpose of this paper is two-fold. First,  we will try to sharpen his
bound with both $d$ and $m$ fixed. Second, we will try to extend his result
to other surfaces.

By dimension count, one may expect that $i(d, m) = dm -
r_{d,m}$ where $r_{d,m}$ is the dimension of the linear series cut out on
$D$ by all curves of degree $m$, namely, $r_{d,m} = m(m+3)/2$ for $m < d$
and $r_{d, m} = dm - (d-1)(d-2)/2$ for $m \ge d$. However, this is simply false
for $m \ge d \ge 3$ by the following construction.

Let $L$ be a bitangent (or flex) line to $D$. Since $D$ is general, $L$ meets
$D$ at $d-2$ distinct points. Let $L(X, Y, Z)$ and $D(X, Y, Z)$ be
the homogeneous defining equations of $L$ and $D$, respectively. Then for
any degree $m - d$ homogeneous polynomial $G(X, Y, Z)$, $L^m(X, Y, Z) +
D(X, Y, Z) G(X, Y, Z) = 0$ defines a degree $m$ curve $C$ which meets $D$ at
$d-2$ points, which are the intersections between $L$ and $D$. If we choose
$G(X, Y, Z)$ general enough, $C$ is irreducible and actually smooth. Hence,
by Xu's result, $i(d, m) = d - 2$ for $m \ge d \ge 3$.

Nevertheless, we think that $i(d,m)$ has the expected value for $d > m$, i.e.,

\begin{conj}\label{conj1}
For $d > m$ and $d\ge 3$, $i(d, m) = dm - r_{d,m}$.
\end{conj}

Although we cannot prove the above conjecture, we have the following
estimate for $i(d,m)$ when $m < d$.

\begin{thm}\label{t1}
For $d > m$,
$$i(d, m) \ge \min\left(dm - \frac{m(m+3)}{2}, 2dm - 2m^2 -2\right).$$
\end{thm}

An easy corollary of the above theorem is the following

\begin{cor}
For $2d \ge 3m - 2$ and $d\ge 3$, 
$i(d, m) = dm - m(m+3)/2$, i.e., \conjref{conj1} holds
for $2d \ge 3m -2$. In particular, it holds for $m \le 4$.
\end{cor}

In order to formulate Kobayashi type conjectures on surfaces other than
$\P^2$, we need to study \questionref{q1} in the following general setting.

\begin{question}\label{q2}
Let $S$ be a smooth surface and let $L$ and $M$ be two line bundles on
$S$. Let $D$ be a general member of $|L|$. What is the minimal number
$i(L, M)$ of points in the set-theoretical intersection $C\cap D$ for any
irreducible curve $C\in |M|$ (suppose that $C$ and $D$ meet properly)?
\end{question}

We will work on rational ruled surfaces, although our method
can be extended to other surfaces. By convention, let $\F_n$ be the
rational ruled surface given by $\P(\CO_{\P^1}\oplus \CO_{\P^1}(n))$
over $\P^1$ and let $C$ and $F$ be the zero section and the fiber of
$\F_n\to \P^1$, i.e., $C^2 = -n$, $C\cdot F = 1$ and $F^2 =
0$. We have the following lower bound for $i(L, M)$ with $L$ ample.

\begin{thm}\label{t2}
Let $L = \CO(aC + bF)$ be an ample line bundle on $\F_n$ with $a\ge 2$ and
$b\ge 2$. Then
$\min_M i(L, M) = \min(a - 1, b - an, b - n - 1)$, where $M$ runs over all
line bundles with irreducible general global sections.
\end{thm}

It follows immediately from \thmref{t2} that every curve on $\F_n$ meets
$D$ at no less than three distinct points for a sufficiently general $D\in
|aC + bF|$ with $a\ge 4$ and $b\ge \max(4+n, 3+an)$. Therefore

\begin{cor}
For a sufficiently general curve $D\in |aC + bF|$ on $\F_n$ with $a\ge 4$ and
$b\ge \max(4+n, 3+an)$, the complement $\F_n\backslash D$ is algebraic
hyperbolic.
\end{cor}

Notice that the bound in \thmref{t2} can be achieved by a curve in $|C|$
or $|F|$, which is necessarily a rational curve. So the lower bounds for
$a$ and $b$ in the above corollary cannot be improved.

This enables us to formulate Kobayashi conjecture on $\F_n$.

\begin{conj}[Kobayashi Conjecture on Rational Ruled Surfaces]
For a sufficiently general curve $D\in |aC + bF|$ on $\F_n$ with $a\ge 4$ and
$b\ge \max(4+n, 3+an)$, the complement $\F_n\backslash D$ is hyperbolic.
\end{conj}

The organization of this paper is as follows. \thmref{t1} and \ref{t2} will
be proved in the next two sections, respectively. At the end of the third
section, we will also discuss some related problems.

Throughout this paper we work exclusively over the field of complex numbers
$\BC$.

\medskip\noindent{\bf Acknowledgments.}
I am very grateful to M. Green, who introduced me to the subject of
hyperbolic geometry and Kobayashi Conjecture, and to C. Hacon for helpful
conversations.

\section{Proof of \thmref{t1}}

Let $W_\delta \subset |\CO(m)|\times |\CO(d)|$ be the incidence correspondence
defined by 
\[
\begin{split}
W_\delta = \{ &(C, D): C\in |\CO(m)| \text { is irreducible}, D\in |\CO(d)|
\text{ is smooth},\\ 
& \text{and $C$ and $D$ meet set-theoretically at $\delta$ points} \}.
\end{split}
\]

Our proof of \thmref{t1} is carried out
by estimating the dimension of $W_\delta$ and
show that $\dim W_\delta < \dim |\CO(d)|$ if $\delta < dm - m(m+3)/2$ and
$\delta < 2dm - 2m^2 - 2$ and hence it cannot dominate $|\CO(d)|$ in this
case.

Let $\pi : W_\delta \to |\CO(m)|$ be the projection of $W_\delta$ to
$|\CO(m)|$ and let $C$ be a general point of $\pi(W_\delta)$ (by a general
point, we mean a general point of an irreducible component of
$\pi(W_\delta)$). 

Let $\pi_C$ be the fiber of $\pi: W_\delta \to |\CO(m)|$ over $C$ and let
$(C, D)$ be a general point on $\pi(C)$.
There exists a series of blowups of $\P^2$
such that the proper transforms $\wt{C}$ and $\wt{D}$ of $C$ and $D$ meet
at smooth points on both curves (since we assume that $D$ is smooth, we
only have to resolve the singularities of $C$ where $D$ passes through).
Let $\wt{\P^2}$ be the resulting blowup of $\P^2$ and $E_i$ ($1\le i \le
\alpha)$ be the
exceptional divisors. Suppose that $\wt{C} \in |\CO(mH - \sum_{i=1}^\alpha 
r_i E_i)$
and $\wt{D}\in |\CO(d H - \sum_{i=1}^\alpha E_i)|$, 
where $r_i > 1$ and $H$ is the
pull-back of the hyperplane divisor of $\P^2$.

Suppose that $\wt{C}$ and $\wt{D}$ meet at points $p_1, p_2, ..., p_\beta$
with multiplicities $m_1, m_2, ...,m_\beta$, respectively, where $\beta
\le \delta\le \alpha + \beta$. 
Then by a deformation-theoretical argument,
the tangent space $T_{\pi_C, (C, D)}$ of $\pi_C$ at $(C, D)$ is contained
in
\[
\begin{split}
&\quad H^0(\CO_{\wt{D}}(dH - \sum_{i=1}^\alpha E_i)\tensor 
\CO_{\wt{D}}(-\sum_{j=1}^\beta (m_j -1) p_j)) \\
&= H^0(\CO_{\wt{D}} ((d-m) H + 
\sum_{i=1}^\alpha (r_i - 1) E_i) \tensor \CO_{\wt{D}}(\sum_{j=1}^\beta p_j)).
\end{split}
\]

Hence by Riemann-Roch

\[
\begin{split}
\dim \pi_C &\le h^0(\CO_{\wt{D}} ((d-m) H + 
\sum_{i=1}^\alpha (r_i - 1) E_i + \sum_{j=1}^\beta p_j))\\
&= \frac{d(d+3)}{2} + \sum_{i=1}^\alpha (r_i -1) + \beta - dm\\
& + h^0(\CO_\wt{D}((m-3)H - \sum_{i=1}^\alpha (r_i - 1) E_i) \tensor 
\CO_{\wt{D}} (-\sum_{j=1}^\beta p_j)).
\end{split}
\]

It is not hard to see that
\begin{equation}\label{e1}
\begin{split}
&\quad h^0(\CO_\wt{D}((m-3)H - \sum_{i=1}^\alpha (r_i - 1) E_i) \tensor 
\CO_{\wt{D}} (-\sum_{j=1}^\beta p_j))\\
&\le h^0(\CO_\wt{C}((m-3)H - \sum_{i=1}^\alpha (r_i - 1) E_i) \tensor 
\CO_{\wt{C}} (-\sum_{j=1}^\beta p_j))\\
&\quad + \sum_{i=1}^\alpha \frac{(r_i-1)(r_i -2)}{2}.
\end{split}
\end{equation}

This can be shown by the following argument.

We further blow up $\wt{\P^2}$ at points $p_1, p_2, ..., p_\beta$ with
corresponding exceptional divisors $F_1, F_2, ..., F_\beta$. We still denote
the resulting surface by $\wt{\P^2}$ and the proper transforms of $C$
and $D$ by $\wt{C}$ and $\wt{D}$. We have the following
exact sequence on $\wt{\P^2}$

\[
\begin{split}
0 &\to H^0((m-d-3)H - \sum_{i=1}^\alpha (r_i-2) E_i)\\
&\to 
H^0((m-3)H - \sum_{i=1}^\alpha (r_i - 1) E_i - \sum_{j=1}^\beta
F_j)\\
&\to H^0(\CO_\wt{D}((m-3)H - \sum_{i=1}^\alpha (r_i - 1) E_i -
\sum_{j=1}^\beta F_j))\\
&\to H^1((m-d-3)H - \sum_{i=1}^\alpha (r_i-2) E_i).
\end{split}
\]

Obviously, $h^0((m-d-3)H - \sum_{i=1}^\alpha (r_i-2) E_i) = 0$ and
\[
\begin{split}
&\quad h^2((m-d-3)H - \sum_{i=1}^\alpha (r_i-2) E_i)\\
&= h^0((d-m)H + \sum_{i=1}^\alpha (r_i -1) E_i + \sum_{j=1}^\beta F_j)\\
&= \frac{(d-m)(d-m+3)}{2} + 1.
\end{split}
\]
Hence by Riemann-Roch, 
$$h^1((m-d-3) H - \sum_{i=1}^\alpha (r_i-2) E_i) =
\sum_{i=1}^\alpha \frac{(r_i-1)(r_i -2)}{2}.$$
Therefore
\begin{equation}\label{e2}
\begin{split}
&\quad h^0(\CO_\wt{D}((m-3)H - \sum_{i=1}^\alpha (r_i - 1) E_i -
\sum_{j=1}^\beta F_j))\\
& \le 
h^0((m-3)H - \sum_{i=1}^\alpha (r_i - 1) E_i - \sum_{j=1}^\beta F_j)
+ \sum_{i=1}^\alpha \frac{(r_i-1)(r_i -2)}{2}.
\end{split}
\end{equation}

Similarly, we have
\begin{equation}\label{e3}
\begin{split}
&\quad h^0(\CO_\wt{C}((m-3)H - \sum_{i=1}^\alpha (r_i - 1) E_i -
\sum_{j=1}^\beta F_j))\\
&= h^0((m-3)H - \sum_{i=1}^\alpha (r_i - 1) E_i - \sum_{j=1}^\beta F_j).
\end{split}
\end{equation}

Combining \eqref{e2} and \eqref{e3}, we obtain \eqref{e1}. Therefore
\[
\begin{split}
\dim \pi_C &\le \frac{d(d+3)}{2} + \sum_{i=1}^\alpha \frac{r_i(r_i -1)}{2} 
+ \beta - dm\\
& + h^0(\CO_\wt{C}((m-3)H - \sum_{i=1}^\alpha (r_i - 1) E_i) \tensor 
\CO_{\wt{C}} (-\sum_{j=1}^\beta p_j))\\
&= \frac{d(d+3)}{2} + \sum_{i=1}^\alpha \frac{r_i(r_i -1)}{2} 
+ \beta - dm + h^0(\omega_{\wt{C}}\tensor \CO_{\wt{C}} (-\sum_{j=1}^\beta p_j))
\end{split}
\]
where $\omega_{\wt{C}}$ is the dualizing sheaf of $\wt{C}$. 

By
Clifford's theorem (see for example \cite[pp. 107-8]{ACGH}), we have either
\begin{equation*}
h^0(\omega_{\wt{C}}\tensor \CO_{\wt{C}} (-\sum_{j=1}^\beta p_j)) = 0
\end{equation*}
or
\begin{equation*}
h^0(\CO_{\wt{C}} (\sum_{j=1}^\beta p_j))\le \beta/2 + 1.
\end{equation*}
Hence correspondingly, we have either
\[
\dim \pi_C\le \frac{d(d+3)}{2} + \sum_{i=1}^\alpha \frac{r_i(r_i -1)}{2} 
+ \beta - dm
\]
or
\[
\dim \pi_C\le \frac{d(d+3)}{2} + \frac{\beta}{2} - dm +
\frac{(m-1)(m-2)}{2}.
\]

Since $C$ is a general member of $\pi(W_\delta)$ and $C$ has singularities
with multiplicities $r_i$ ($1\le i \le \alpha$), by Zariski's theorem on the
deformation of planary curve singularities \cite{Z}, we have
$$\dim \pi(W_\delta) \le \frac{m(m+3)}{2} - \sum_{i=1}^\alpha
\frac{r_i(r_i -1)}{2}.$$
And hence we have either
\[
\dim W_\delta \le \frac{d(d+3)}{2} + \beta - dm + \frac{m(m+3)}{2}
\]
or
\[
\dim W_\delta \le \frac{d(d+3)}{2} + \frac{\beta}{2} -dm +
m^2 + 1 - \sum_{i=1}^\alpha \frac{r_i(r_i -1)}{2}.
\]
Therefore, if $W_\delta$ dominates $|\CO(d)|$, we necessarily have
$$\delta \ge \beta \ge \min\left(dm - \frac{m(m+1)}{2}, 
2(dm - m^2 -1)\right).$$
This finishes the proof of \thmref{t1}.

\section{Intersections of Two Curves on Rational Ruled Surfaces}

Our approach to \thmref{t2} is different from that of Xu's. A key
ingredient of Xu's proof of \thmref{txu} is a map from the deformation
space of the pair $(D, E)$, where $D\in |\CO_{\P^2}(d)|$ and $E\in
|\CO_{\P^2}(m)|$ meet at no less than $s$ distinct points, to the
cohomology group of a sheaf over $D$. More specifically, let $(Z_0, Z_1,
Z_2)$ be generic homogeneous
coordinates of $\P^2$ and let $F_0\in H^0(\CO_{\P^2}(d))$ and $G_0\in 
H^0(\CO_{\P^2}(m))$ be the defining equations of $D$ and $E$. A first order
deformation of $(D, E)$ is given by $F_1\in H^0(\CO_{\P^2}(d))$ and $G_1\in
H^0(\CO_{\P^2}(m))$ such that the curves $\{F_0 + tF_1 = 0\}$ and $\{G_0 +
tG_1 = 0\}$ meet at no less than $s$ points over the ring
$\BC[t]/(t^2)$. It is observed by Xu that \cite[Lemma 1]{X}
\begin{equation}\label{e4}
\frac{\partial F_0}{\partial Z_i} G_1 - \frac{\partial G_0}{\partial Z_i}
F_1 \in H^0(D, \CO_D(d + m - 1)\tensor \CO_D(-\sum_{j=1}^s (\mu_j - 1)p_j))
\end{equation}
for $i = 0, 1, 2$, where $D$ and $E$ meet at $p_1, p_2, ..., p_s$ with
multiplicities $\mu_1, \mu_2, ...,\mu_s$, respectively.

The relation \eqref{e4} forms the basis of Xu's proof of \thmref{txu}. If
we were to prove \thmref{t2} following Xu's line of argument, we would have
to come up with a relation similar to \eqref{e4} on $\F_n$, which we are
unable to do. So we find that Xu's analysis, though ingenious on its own,
is hard, if not impossible, to carry out on surfaces other than
$\P^2$. Therefore, we will adopt a different approach to \thmref{t2}, which
is based upon degeneration and induction.

Let $\Delta$ be a disk parameterized by $t$ and let $Y\subset \F_n\times \Delta$
be a pencil of curves in $|aC + bF|$ whose central fiber $Y_0 = G\cup \Gamma$
is reducible with two components $G\in |(a-1)C + (b-n-1)F|$ and $\Gamma
\isom \P^1\in |C + (n+1)F|$. Let $X\subset \F_n\times \Delta$ be a family
of curves on $\F_n$ whose general fiber $X_t$ ($t\ne 0$) meets $Y_t$ at $s$
distinct points (a base change may be needed to ensure the existence of
$X$). If $X_0$ meets $Y_0$ properly, we may deduce $s \ge \min(a - 1, b -
an, b - n - 1)$ by the induction hypothesis that
\[
\begin{split}
\# (X_0\cap G) &\ge \min(a - 2, b - an - 1, b - n - 2)\\
&= \min(a - 1, b - an, b - n - 1) - 1
\end{split}
\]
and by the fact that $\# (X_0\cap \Gamma)\ge 1$, where we use the notation
$\# (A\cap B)$ to denote the number of points in the set-theoretical
intersection $A\cap B$ between the two curves $A$ and $B$.
Of course, some care has
to be taken to make sure that $X_0$ meets $\Gamma$ at at least one point
outside of $G\cap \Gamma$ (see below). Unfortunately, $X_0$ may very well
contain $G$ or $\Gamma$ as a component. So we have to regard
$|\CO_{Y_0}(X_0)|$ as the limit linear series $\lim_{t\to 0}
|\CO_{Y_t}(X_t)|$ and, correspondingly, $Y_0\cap X_0$ as the limit of the
section $Y_t\cap X_t$ in $|\CO_{Y_t}(X_t)|$. For an introduction to the
theory of limit linear series, please see, for example, \cite{E-H} or
\cite[Chap. 5]{H}.

For the purpose of induction, we will prove \thmref{t2} in the following
slightly more general form.

\begin{prop}\label{prop1}
Let $L = \CO(aC + bF)$ be an ample line bundle on $\F_n$ with $a\ge 2$ and
$b\ge 2$. Then for a sufficiently general curve $D\in |L|$,
\begin{enumerate}
\item
$\# (D\cap E)\ge
\min(a - 1, b - an, b - n - 1)$ for any curve $E\subset \F_n$ that meets
$D$ properly;
\item
in addition, there exists a set $\Sigma_D$ consisting of
countably many points on $D$ such that if $\# (D\cap E) = \min(a - 1, b -
an, b - n - 1)$ for some $E$, $(D\cap E) \subset \Sigma$.
\end{enumerate}
\end{prop}

We prove \propref{prop1} by induction on $\min(a - 1, b - an, b - n - 1)$.

If $\min(a - 1, b - an, b - n - 1) = 1$, we only need to verify the second
part of the proposition. Notice that $D$ has genus $g(D) = 1 + \frac{1}{2}(a-2)(b -
an) + \frac{1}{2} a(b - n -2)\ge 1$. If $D$ meets $E$ at a single $p$ for
some $E$, $\CO_D(\mu p) = \CO_D(E)$, where $\mu = D\cdot E$. If we fix the
divisor class of $E$, there are only finitely many points $p$ with this
property since $g(D)\ge 1$. Therefore, there are only countably many points
$p$ such that $D\cap E = \{p\}$ for some $E$.

Suppose that $\min(a - 1, b - an, b - n - 1) \ge 2$. Notice that $\CO((a-1)C
+ (b - n - 1)F)$ is ample under this assumption.

Let $X, Y, G$ and $\Gamma$ be defined as before. Suppose that $G$ and
$\Gamma$ meet at points $p_1, p_2, ..., p_l$, where $l = a + b - n - 2$.
Let $M = \CO(X_t)$ be the line bundle associated to $X_t$.

Let $\sigma_t\in |\CO_{Y_t}(X_t)|$ be the section cut out by $X_t$ on $Y_t$
and let $\sigma_0 = \lim_{t\to 0} \sigma_t$. 
Let $\sigma_\Gamma = \sigma_0 |_\Gamma$ and $\sigma_G = \sigma_0|_G$ be the
restrictions of $\sigma_0$ to $\Gamma$ and $G$, respectively.
Then $\sigma_\Gamma$ is a section in
\[
|\CO_\Gamma(\mu(p_1 + p_2 + ... + p_l))\tensor M| = |\CO_\Gamma(\mu
G)\tensor M|
\]
and $\sigma_G$ is a section in
\[
|\CO_G(-\mu(p_1 + p_2 + ... + p_l))\tensor M| = |\CO_G(-\mu \Gamma)\tensor M|
\]
where $\mu$ is an integer and $\sigma_\Gamma$ and $\sigma_G$
are cut out by sections in $|\CO(\mu G) \tensor M|$ and
$|\CO(-\mu \Gamma)\tensor M|$, respectively.

Suppose that $\CO(-\mu \Gamma)\tensor M$ is nontrivial. Then by induction
hypothesis $\sigma_G$ vanishes at no less than $\min(a - 2, b - an - 1,
b - n - 2)$ distinct points. If $\sigma_\Gamma$ vanishes at at least one
point other than $p_1, p_2, ..., p_l$, we are done; if not, we have
either $\CO(\mu G) \tensor M$ is trivial and $\sigma_\Gamma$ is nowhere
vanishing or $\sigma_\Gamma$ only vanishes at $p_1, p_2, ..., p_l$.

If $\CO(\mu G) \tensor M$ is trivial and $\sigma_\Gamma$ is nowhere
vanishing, then for any two points among $p_1, p_2, ..., p_l$, say $p_1$
and $p_2$, the ratio $\sigma_G(p_1)/\sigma_G(p_2)$ is
uniquely determined by the choice of the pencil $Y$. Actually we have the
following very explicit relation
\begin{equation}\label{e5}
\frac{\sigma_G(p_1)}{\sigma_G(p_2)}
= \left(\frac{f(p_1)}{f(p_2)}\right)^{-\mu}
\end{equation}
where $f\in |L|$ is the section which cuts out a general member $Y_t$ of the
pencil $Y$.
If $\sigma_G$ vanishes at more than $\min(a - 2, b - an - 1,
b - n - 2)$ distinct points, there is nothing to prove;
otherwise, $\sigma_G$ vanishes at exactly $\min(a - 2, b - an - 1,
b - n - 2)$ distinct points. Then by induction hypothesis, there are only
countably many possible choices of $\sigma_G$. However, by
\eqref{e5}, the ratio $\sigma_G(p_1)/\sigma_G(p_2)$ can be made into an
arbitrary complex value by a choice of $f$ (and thus a choice of the pencil
$Y$). Contradiction.

If $\sigma_\Gamma$ only vanishes at $p_1, p_2, ..., p_l$, since we have
already taken care of the case that $\CO(\mu G) \tensor M$ is trivial and
$\sigma_\Gamma$ is nowhere vanishing, we may assume that $\sigma_\Gamma$
vanishes at at least one point among $p_1, p_2, ..., p_l$, say $p_1$.
Then $\sigma_G$ must vanish at $p_1$ as well. Again, if $\sigma_G$ vanishes
at more than $\min(a - 2, b - an - 1, b - n - 2)$ distinct points, there is
nothing to prove; otherwise, $\sigma_G$ vanishes at exactly $\min(a - 2, b
- an - 1, b - n - 2)$ distinct points. By induction hypothesis, $p_1\in
\Sigma_G$. But if we choose $\Gamma$ generically, $p_1\not\in
\Sigma_G$. Contradiction.

Now suppose that $\CO(-\mu \Gamma)\tensor M$ is trivial. If $\sigma_G = 0$,
then $\sigma_\Gamma$ vanishes at $p_1, p_2, ..., p_l$ and $l = a + b - n -
2 > \min(a - 1, b - an, b - n - 1)$; we are done. Otherwise, $\sigma_G$ is
no where vanishing. The ratio $\sigma_\Gamma(p_i)/
\sigma_\Gamma(p_j)$ for any two points $p_i$ and $p_j$ among $p_1, p_2,
..., p_l$, just as in \eqref{e5}, is uniquely determined by the choice of
$Y$ and is given by
\begin{equation}\label{e6}
\frac{\sigma_\Gamma(p_i)}{\sigma_\Gamma(p_j)}
= \left(\frac{f(p_i)}{f(p_j)}\right)^{\mu}.
\end{equation}
The rational map $|L|\to \P^{l-1}$ by sending $f\in |L|$ to
\begin{equation}\label{e7}
\left(f^\mu(p_1), f^\mu(p_2), ..., f^\mu(p_l)\right)
\end{equation}
is dominant due to the facts that $H^0(\F_n, L)$ surjects
onto $H^0(\Gamma, L)$ and 
$L\tensor \CO_\Gamma(-\sum_{i\ne j} p_i)$
is base point free on $\Gamma$ for each $1\le j \le l$.
On the other hand, the space
\[
\{\sigma_\Gamma | \sigma_\Gamma \text{ vanishes at less than } l - 1 \text{
distinct points}\}
\]
has dimension $l-2$ and hence cannot dominate $\P^{l-1}$. So
$\sigma_\Gamma$ vanishes at at least $l - 1 > \min(a - 1, b - an, b - n -
1)$ distinct points for a general choice of $f$ by \eqref{e6}.

This finishes the proof of the first part of the proposition.

Suppose that $\sigma_0$ vanishes at exactly $\min(a - 1, b - an, b - n -
1)$ distinct points. This can happen only when $\CO(-\mu \Gamma)\tensor M$
is nontrivial.

Suppose that $\sigma_G$ vanishes at exactly $\min(a - 2, b - an-1, b - n -
2)$ distinct points. Our previous argument shows that $\sigma_G$ does not
vanish at $p_1, p_2, ..., p_l$ for a general choice of $G\cup \Gamma$.
Then $\sigma_\Gamma$ must vanish at a single point $p\not\in \{p_1, p_2,
..., p_l\}$. Since $\# (G\cap \Gamma) \ge 2$, the natural
map from $Y_0\backslash \{p_1,
p_2, ..., p_l\}$ to $\Pic(Y_0)$ is injective. So $p$ is determined up to
finitely many possibilities by $M$ and the vanishing locus of
$\sigma_G$. By induction, the vanishing locus of $\sigma_G$ is contained in
some countable set $\Sigma_G$ depending only on $G$.
So the vanishing locus of $\sigma_0$ is also
contained in some countable set $\Sigma_{G\cup\Gamma}$ depending only on
$G\cup\Gamma$.

Suppose that $\sigma_G$ vanishes at exactly $\min(a - 1, b - an, b - n -
1)$ distinct points and suppose that there is a one-parameter family of
$\sigma_0(u)$ with this property, where $\sigma_0(u)$ is parameterized by
$u\in U$ for some irreducible curve $U$.

Suppose that $\CO(\mu G) \tensor M$ is trivial. There exists $u_0\in U$
such that $\sigma_G(u_0)$ vanishes at $p_1$. Since $\CO(\mu G) \tensor M$
is trivial, $\sigma_\Gamma(u_0) = 0$ and hence $\sigma_G(u_0)$ vanishes at
$p_1, p_2, ..., p_l$. But $l > \min(a - 1, b - an, b - n - 1)$. Contradiction.

Suppose that $\CO(\mu G) \tensor M$ is nontrivial. Then $\sigma_\Gamma(u)$
vanishes at at least one point among $p_1, p_2, ..., p_l$, say $p_1$. Hence
$\sigma_G(u)$ vanishes at $p_1$ for all $u\in U$. As $u$ varies, another
vanishing point of $\sigma_G(u)$ will approach $p_1$. So there exists $u_0\in
U$ such that $\sigma_G(u_0)$ vanishes at $\min(a - 2,
b - an - 1, b - n - 2)$ distinct points and among them vanishes at a
general point $p_1$. Again this is impossible by induction. Contradiction.

This finishes the proof of \propref{prop1}.

The degeneration method we used can be applied to surfaces other than
rational ruled surfaces. For example, we can give an alternative proof of Xu's
\thmref{txu} by degenerating a degree $d$ curve to a union of a degree
$d-1$ curve and a line and arguing by induction.

\begin{proof}[A proof of Xu's \thmref{txu} via degeneration]
As in the case of \propref{prop1}, we need to add a clause to the theorem
for the purpose of induction, i.e., we will prove the following statement by
induction on $d$.

For a sufficiently general curve $D$ of degree $d \ge 3$ in $\P^2$,
$\# (D\cap E) \ge d - 2$ for any curve $E\subset \P^2$ that meets $D$
properly. In addition, there exists a set $\Sigma_D$ of countably many
points on $D$ such that if $\# (D\cap E) = d - 2$ for some $E$, $(D\cap E)
\subset \Sigma_D$.

Let $Y\subset \P^2\times \Delta$ be a pencil of degree $d$ curves whose
central fiber $Y_0 = G\cup \Gamma$ is the union of a curve $G$ of degree
$d-1$ and a line $\Gamma$ and let $G\cap \Gamma = \{ p_1, p_2, ..., p_l\}$
where $l = d - 1$.

Let $X, M, \sigma_t, \sigma_0, \sigma_G, \sigma_\Gamma$ and $\mu$ be defined as
before. Almost nothing in the argument of \propref{prop1} needs
changing except in the case that $\CO(-\mu\Gamma)\tensor M$ is trivial and
$\sigma_G$ is nowhere vanishing. In this case, following our previous
argument, we can show that $\sigma_\Gamma$ vanishes at no less than $l - 1$
points. The difference is that now we have $l - 1 = d - 2$ and we have to
verify that there are only finitely many $\sigma_\Gamma$ that vanishes at exactly
$l - 1$ distinct points. This is more or less obvious because the map from
$|L|$ to $\P^{l-1}$ given by \eqref{e7} is dominant and the space
\[
\{\sigma_\Gamma | \sigma_\Gamma \text{ vanishes at exactly } l - 1 \text{
distinct points}\}
\]
has dimension $l-1$.
\end{proof}

Our degeneration method also works for Del Pezzo surfaces.

\begin{thm}\label{t3}
Let $\wt{\P^2}$ be the blowup of $\P^2$ at $2\le r \le 6$ general points and let
$L_1, L_2, ..., L_k,...$ be all the smooth rational curves on $\wt{\P^2}$ with
self-intersection $-1$. Let $L$ be an ample line bundle on
$\wt{\P^2}$. Then for a sufficiently general curve $D\in |L|$,
\begin{enumerate}
\item $\# (D\cap E) \ge \min_{k} (D\cdot L_k)$ for any curve $E\subset
\wt{\P^2}$ that meets $D$ properly;
\item in addition, there exists a set $\Sigma_D$ of countably many
points on $D$ such that if $\# (D\cap E) = \min_{k}(D\cdot L_k)$ for some
$E$, $(D\cap E) \subset \Sigma_D$.
\end{enumerate}
Therefore, for a sufficiently general curve $D\in |L|$ with
$\min_{k} (D\cdot L_k) \ge 3$, the complement $\wt{\P^2}\backslash D$ is
algebraic hyperbolic.
\end{thm}

\begin{proof}
Let $K_{\wt{\P^2}}$ be the canonical divisor of $\wt{\P^2}$.
We argue by induction on $\min_{k} (D\cdot L_k)$.

For $\min_{k} (D\cdot L_k) = 1$, we need to verify that $g(D) \ge 1$, which
is more or less obvious.

Suppose that $\min_{k} (D\cdot L_k) \ge 2$. Let $Y\subset \P^2\times \Delta$
be a pencil of curves in $|L|$ whose central fiber $Y_0 = G\cup\Gamma$ is a
union of $G\in \left|L\tensor \CO(K_{\wt{\P^2}})\right|$ and $\Gamma\in
\left|-K_{\wt{\P^2}}\right|$ and let $G\cap \Gamma = \{p_1, p_2, ..., p_l\}$.

Let $X, M, \sigma_t, \sigma_0, \sigma_G, \sigma_\Gamma$ and $\mu$ be defined as
before. Again, the same argument for \propref{prop1} goes through.
We need only to check the following facts, all of which are routine exercises.
\begin{enumerate}
\item $l > \min_{k} (D\cdot L_k)$.
\item $H^0(\wt{\P^2}, L)$ surjects onto $H^0(\Gamma, L)$ and 
\[
L\tensor
\CO_\Gamma(-\sum_{i\ne j} p_i) = \CO_\Gamma(-K_{\wt{\P^2}})\tensor
\CO_\Gamma(p_j)
\]
is  base point free on $\Gamma$ for each $1\le j\le l$.
Hence the map from $|L|$ to $\P^{l-1}$ given by
\eqref{e7} is dominant.
\item In the case that $\CO(-\mu\Gamma)\tensor M$ is trivial and
$\sigma_G$ is nowhere vanishing, we can prove that $\sigma_\Gamma$ vanishes
at no less than $l-1$ distinct points as before. But actually, we can do
better here since the space
\[
\{\sigma_\Gamma | \sigma_\Gamma \text{ vanishes at less than } l \text{
distinct points}\}
\]
has dimension $l-2$ due to the fact that $\Gamma$ is elliptic instead of
rational. Therefore, $\sigma_\Gamma$ vanishes at no less than $l$ distinct
points.
\end{enumerate}
\end{proof}

When we go up in dimension, however, some essential difficulties
present themselves. For example, in $\P^3$, fix a sufficient general
surface $S$ of degree $d$ and it is expected that any curve meets $S$ at no
less than $d - 4$ distinct points \cite[Question 2]{X}.
Let $Y$ be a pencil of degree $d$ surfaces whose
central fiber is a union of a degree $d-1$ surface and a plane and
let $X$ be family of curves in $\P^3$ meeting $Y$ fiberwise. To carry out
the argument as in dimension two, we need to take the limit $X_t\cap Y_t$
as an element in $A_0(Y_t)$, the 0-dimension Chow ring of $Y_t$. Of course,
we do not know how to do this at present.

\end{document}